\documentclass[12pt]{amsart}

\usepackage{latexsym}
\usepackage[latin1]{inputenc}
\usepackage{amssymb,amsmath,amsfonts}
\usepackage{enumerate}

\newtheorem{MainTheorem}{Theorem}

\newtheorem{Proposition}{Proposition}[section]
\newtheorem{Definition}[Proposition]{Definition}
\newtheorem{Lemma}[Proposition]{Lemma}
\newtheorem{Theorem}[Proposition]{Theorem}

\DeclareMathOperator{\vol}{Vol}
\DeclareMathOperator{\kl}{Kl}
\DeclareMathOperator{\Val}{Val}

\DeclareMathOperator{\Kl}{Kl}
\DeclareMathOperator{\Gr}{Gr}
\DeclareMathOperator{\spa}{span}
\DeclareMathOperator{\SL}{SL}
\DeclareMathOperator{\GL}{GL}
\DeclareMathOperator{\SU}{SU}

\newcommand{\R}{\mathbb{R}}
\newcommand{\C}{\mathbb{C}}

\title{Projection bodies in complex vector spaces}
\author{Judit Abardia} 
\author{Andreas Bernig}
\address{Institut f\"ur Mathematik, Goethe-Universit\"at Frankfurt, 
Robert-Mayer-Str. 10, 60054 Frankfurt, Germany}

\email{abardia@math.uni-frankfurt.de}
\email{bernig@math.uni-frankfurt.de}

\begin{document}
\begin{abstract}
The space of Minkowski valuations on an $m$-dimensional complex vector space which are continuous, translation
invariant and contravariant under the complex special linear group is explicitly described.
Each valuation with these properties is
shown to satisfy geometric inequalities of Brunn-Minkowski, Aleksandrov-Fenchel and Minkowski type. 
\end{abstract}

\thanks{{\it MSC classification}:  52B45, 52A39, 52A40}

\maketitle 

\section{Introduction}

Let $V$ be a real vector space of dimension $n$. Let $\mathcal{K}(V)$ denote the space of
non-empty compact convex bodies in $V$, endowed with the Hausdorff topology. 

The projection body of $K \in \mathcal{K}(V)$ is the convex body $\Pi K \in \mathcal{K}(V^*)$ whose support function is given by 
\begin{displaymath}
h(\Pi K,v)=\frac{n}{2} V(K,\ldots,K,[-v,v]), \quad \forall v \in V.  
\end{displaymath}
Here $V(K,\ldots,K,[-v,v])$ is the mixed volume of $n-1$ copies of $K$ and one copy of the segment $[-v,v]$ joining $-v$ and $v$. In the following, we will use the standard notation $V(K[n-1],[-v,v])$ instead of $V(K,\ldots,K,[-v,v])$.

In more intuitive terms, suppose that $V$ is endowed with a Euclidean scalar product. Then we can identify
$V^*$ with $V$ and the support function of $\Pi K$ in the direction $v \in S^{n-1}$ is the volume of the
orthogonal projection of $K$ onto the hyperplane $v^\perp$. 

Projection bodies have been widely studied since their introduction by Minkowski at the end of 19th century. 
They satisfy important properties which make them useful not only in convex geometry, but also in other areas such as geometric tomography, stereology, computational geometry, optimization or functional analysis (see, for example,  \cite{bolker,bourgain_lindenstrauss,gardner_book06, koldobsky, lutwak93,schneider82,stanley}). 

As an example, we mention the solution of Shephard's problem \cite{goodey_zhang96,petty67,schneider67,shephard64}, where projection bodies played an important role. 

There are also important inequalities involving the volume of the projection body and its polar, such as the Petty projection inequality \cite{petty71} and the Zhang projection inequality \cite{zhang91}. For additional information and recent results on projection bodies see, for example, \cite{gardner_book06,haberl_schuster09b,leichtweiss_book98,lutwak93,lutwak_handbook93,lutwak_yang_zhang00,lutwak_yang_zhang10,schneider_book93,thompson_book96}.

Ludwig \cite{ludwig02, ludwig_2005} proved that the projection operator $\Pi$, sending each convex body to its projection body is characterized by the following properties: 
\begin{enumerate}
 \item $\Pi$ is a continuous Minkowski valuation, i.e. 
\begin{displaymath}
 \Pi(K \cup L)+\Pi(K \cap L)=\Pi K + \Pi L
\end{displaymath}
whenever $K, L, K \cup L \in \mathcal{K}(V)$ (here the sum is the Minkowski sum of convex bodies);
\item $\Pi$ is translation invariant, i.e. $\Pi(K+x)=\Pi(K)$ for all $x \in V$;
\item $\Pi$ is $\SL(V)$-contravariant, i.e. 
\begin{displaymath}
 \Pi(gK)=g^{-*}\Pi(K),\quad \forall K\in \mathcal{ K}(V), g\in \SL(V).
\end{displaymath}
\end{enumerate}

More generally, a valuation is an operator $Z:\mathcal{K}(V) \to (A,+)$ with $(A,+)$ an abelian semi-group, such that the following additivity property is satisfied: 
\begin{displaymath}
Z K_1+ZK_2=Z(K_1\cup K_2)+Z(K_1\cap K_2),\quad \forall K_1,K_2,K_1\cup K_2\in\mathcal{K}(V).
\end{displaymath}

The classical case $A=\mathbb{R}$ has attracted a lot of attention, we refer to \cite{alesker_fourier,alesker_mcullenconj01, alesker03_un, alesker04_product,bernig_broecker07,bernig_fu_hig,fu06,klain00,ludwig_reitzner99} for some new developments. 

In the case of the projection body, $A=\mathcal{K}(V^*)$ endowed with Minkowski addition. Valuations with values
in $\mathcal{K}(V)$ or $\mathcal{K}(V^*)$ are called {\it Minkowski valuations}. See
\cite{alesker_bernig_schuster,haberl10,ludwig10,ludwig02,ludwig_2005,ludwig06,ludwig06_survey,schneider_schuster06,schuster08, schuster10, schuster_wannerer10} for more
information on Minkowski valuations. 

In this paper we study a complex version of Ludwig's characterization theorem of the projection operator. This
work is part of a larger program aiming at complex-affine versions of some geometric and functional
inequalities (e.g. isoperimetric inequalities and Sobolev inequalities). It seems that not much has been done
in this direction. We know only of one work, namely the solution of a complex version of the
Busemann-Petty problem in \cite{koldobsky_koenig_zymonopoulou08}. 

Let us now describe our main result. The real vector space $V$ of real dimension $n$ is replaced by a complex
vector space $W$ of complex dimension $m$ and the group $\SL(V)=\SL(n,\mathbb{R})$ is
replaced by the group $\SL(W,\mathbb{C})=\SL(m,\mathbb{C})$. Note that $\SL(m,\mathbb{C}) \subset \SL(2m,\mathbb{R})$, so that each element in $\SL(m,\mathbb{C})$ is volume preserving.

For a subset $C \subset \C$ and $w \in W$, we
denote 
\begin{displaymath}
 C \cdot w:=\{cw| c \in C\} \subset W. 
\end{displaymath}
Note that $C\cdot w$ is convex if $C$ is convex. 

\begin{MainTheorem}\label{thm_t1}
Let $W$ be a complex vector space of complex dimension $m \geq 3$. A map $Z:\mathcal{K}(W) \to
\mathcal{K}(W^*)$ is a continuous translation invariant and $\SL(W,\mathbb{C})$-contravariant Minkowski
valuation if and only if there exists a convex body $C\subset\C$ such that $Z=\Pi_C$, where
$\Pi_CK \in \mathcal{K}(W^*)$ is the convex body with support function 
\begin{equation} \label{eq_def_pi_c}
 h(\Pi_{C}K,w)=V(K[2m-1],C\cdot w),\quad\forall w \in W.
\end{equation}
Moreover, $C$ is unique up to translations.
\end{MainTheorem}

Let us point out the complete analogy with the real case: replacing formally $\mathbb{C}$ by $\mathbb{R}$ in
the theorem gives Ludwig's theorem. 

The assumption $m \geq 3$ is essential. In Proposition \ref{prop_n2} we will construct a class of continuous,
translation invariant, $\SL(2,\C)$-contravariant Minkowski valuations on $\C^2$ which are not of the form
\eqref{eq_def_pi_c}.

In the last section of this paper, we establish a number of inequalities for complex projection bodies which are analogs of inequalities satisfied by the classical projection body.

Before stating these theorems, let us define mixed complex projection body. Their real analogous were introduced by Bonnesen-Fenchel \cite{bonnesen_fenchel87} and studied, among others, by Chakerian, Goodey and
Lutwak \cite{chakerian67,goodey86,lutwak85,lutwak86,lutwak93}.

\begin{Definition}\label{MCPB} 
Let $K_1,\dots,K_{2m-1} \in \mathcal{K}(W)$ and $C\subset\C$. The {\em mixed complex projection body
$\Pi_C(K_1,\dots,K_{2m-1}) \in \mathcal{K}(W^*)$} is the convex body whose support function is given by
\begin{displaymath}
 h(\Pi_C(K_1,\dots,K_{2m-1}),w)=V(K_1,\dots,K_{2m-1},C \cdot w), \quad \forall w \in W.
\end{displaymath}
\end{Definition}

We fix a Euclidean scalar product on $W$, denote by $B$ the unit ball and use the following notation:
\begin{displaymath}
 W_i(K,L)=V(K[2m-1-i],B[i],L), \quad K,L \in \mathcal{K}(W), 0 \leq i \leq 2m-1. 
\end{displaymath}

The mixed volume $W_i(K,K)$ will be written as $W_i(K)$ and is called the
{\it $i$-th quermassintegral} of $K$. 

\begin{MainTheorem}\label{thm_geom_ineq}
Let $K,L,K_1,\dots K_{2m-1} \in \mathcal{K}(W)$. 
\begin{enumerate}[i)]
\item Brunn-Minkowski type inequality.\label{thm_BMCT} 
\begin{displaymath}\label{eq_BMC}
\vol(\Pi_C(K+L))^{1/2m(2m-1)} \geq \vol(\Pi_CK)^{1/2m(2m-1)} + \vol(\Pi_CL)^{1/2m(2m-1)}.
\end{displaymath}
\smallskip
\item Aleksandrov-Fenchel type inequality. For $0\leq  i \leq
2m-1, 2\leq k\leq 2m-2$,\label{thm_AFIneq}
\begin{displaymath}
 W_i(\Pi_C(K_1,\dots,K_{2m-1}))^k \geq \prod_{j=1}^k W_i(\Pi_C(K_j,\dots,K_j,K_{k+1},\dots,K_{2m-1})).
\end{displaymath}
\smallskip
\item Minkowski type inequality. For $0 \leq i < 2m-1$, \label{thm_MIneq}
\begin{displaymath}
 W_i(\Pi_C(K[2m-2],L))^{2m-1}\geq W_i(\Pi_C K)^{2m-2}W_i(\Pi_C L).
\end{displaymath}
\end{enumerate}
Moreover, if $K$ and $L$ have non-empty interior and $C$ is not a point then the equality in i) and iii) holds if and only if $K$ and $L$ are homothetic. 
\end{MainTheorem}

\subsection*{Acknowledgments}
We would like to thank Franz Schuster for many interesting discussions and useful remarks on this paper. We also thank the anonymous referee for several helpful suggestions.


\section{Background and conventions}

In this section, we fix some notation which will be used later on. We try to use intrinsic definitions
whenever possible. In particular, we do not assume that our vector space $V$ is endowed with an inner
product, hence we will distinguish between $V$ and $V^*$. 

\subsection{Support function}

The support function of $K \in \mathcal{K}(V)$ is the function on $V^*$ defined by
\begin{align*}
 h_K:V^* & \to \R, \\
\xi & \mapsto \sup_{x \in K}\langle \xi,x\rangle,
\end{align*}
where $\langle \xi,x\rangle$ denotes the pairing of $\xi \in V^*$ and $x\in V$.

The support function is $1$-homogeneous (i.e. $h_K(t\xi)=th_K(\xi)$ for all $t \geq 0$) and subadditive (i.e.
$h_K(\xi+\eta) \leq h_K(\xi)+h_K(\eta)$). Conversely, every $1$-homogeneous and subadditive
function on $V^*$ is the support function of a unique compact convex set $K \in \mathcal{K}(V)$ (cf.
\cite[Theorem 1.7.1]{schneider_book93}). 

Throughout this paper we also shall use the following property of the support function: 
\begin{equation} \label{eq_adjoint}
 h(gK,\xi)=h(K,g^*\xi),\quad \forall \xi \in V^*,g \in \GL(V).
\end{equation}

\subsection{Mixed Volumes}
We refer to \cite{schneider_book93} for details about mixed volumes.   

In a real vector space $V$ of dimension $n$ with a volume $\vol$, the {\em mixed
volume} is the unique symmetric and Minkowski multilinear map $(K_1,\ldots,K_n) \mapsto
V(K_1,\ldots,K_n)$ with $V(K,\ldots,K)=\vol(K)$. 

These functionals are nonnegative, continuous, symmetric and translation invariant in each
component. Moreover,  
\begin{equation*} 
V(gK_1,\dots,gK_n)=|\det g|V(K_1,\dots,K_n), \quad g \in \GL(V),
\end{equation*}
so that they are invariant under volume-preserving affine transformations $g\in SL(V)$.

We will also need the monotonicity property: if $K_1 \subset K_1'$, then 
\begin{equation*} 
V(K_1,\dots,K_n)\leq V(K_1',\dots,K_n).
\end{equation*}

Given convex bodies $K_1,\ldots,K_{n-1}$, there is a unique extension of the functional $K \mapsto
V(K_1,\ldots,K_{n-1},K)$ to a linear functional on the space of all continuous $1$-homogeneous
functions $f:V^* \to \R$ such that for all $L \in \mathcal{K}(V)$ 
\begin{displaymath}
 V(K_1,\ldots,K_{n-1},h_L)=V(K_1,\ldots,K_{n-1},L).
\end{displaymath}

\subsection{Translation invariant valuations}

Let $V$ be a finite-dimensional real vector space of dimension $n$. The Banach space of real-valued, translation
invariant, continuous valuations on $V$ is denoted by $\Val$. It has been studied intensively, see for
instance \cite{alesker_fourier,alesker_mcullenconj01,klain00,mcmullen80}. 

A basic structure result for $\Val$ is the following direct sum decomposition due to McMullen. A valuation $\mu \in \Val$ is
called {\it homogeneous} of degree $k$ if $\mu(tK)=t^k\mu(K)$ for all $t \geq 0$. If $\mu(-K)=\mu(K)$ for all
$K$,
then $\mu$ is called {\it even}; if $\mu(-K)=-\mu(K)$, then $\mu$ is called {\it odd}. The subspace of
even/odd
valuations of degree $k$ is denoted by $\Val_k^\pm$. 

\begin{Theorem}[McMullen \cite{mcmullen77}]
 \begin{equation} \label{eq_mcmullen_dec}
  \Val=\bigoplus_{\substack{k=0,\ldots,n\\ \varepsilon=\pm}} \Val_k^\varepsilon.
 \end{equation} 
\end{Theorem}

The space $\Val_k^+$ admits the following geometric description. In order to simplify the notation, let us fix
a Euclidean scalar product on $V$. Let $\Gr_kV$ be the Grassmannian manifold of all $k$-dimensional
subspaces in $V$. 

Let $\mu \in \Val_k^+$ and let $E$ be a $k$-dimensional subspace of $V$. By a theorem of Klain \cite{klain00}, 
$\mu|_E$ is a multiple of the volume on $E$:
\begin{displaymath}
 \mu(K)=\kl_\mu(E) \vol(K), \quad \forall K \in \mathcal{K}(E).
\end{displaymath}
The function $\kl_\mu:\Gr_k(V) \to \R$ is called {\it Klain function}. Klain's injectivity theorem
\cite[Theorem 3.1]{klain00} states
that the valuation $\mu \in \Val_k^+$ is uniquely determined by its Klain function $\kl_\mu \in C(\Gr_kV)$. 

The next notion which we need comes from representation theory. We refer to \cite{wallach} for more details. 

The group $\GL(V)$ acts naturally on $\Val$:
\begin{displaymath}
 g\mu(K)=\mu(g^{-1}K), \quad g \in \GL(V), K \in \mathcal{K}(V). 
\end{displaymath}

If for a valuation $\mu \in \Val$, the map $g \mapsto g\mu$ from the Lie group $\GL(V)$ to the Banach space
$\Val$ is smooth, then $\mu$ is called {\it smooth}. The subspace of smooth valuations is denoted by
$\Val^{sm}$, it is a dense subspace in $\Val$. 

Clearly, if $\mu \in \Val_k^{sm,+}$, then the Klain function of $\mu$ is smooth as a function on $\Gr_kV$. We
will need the following fact about smooth valuations.

\begin{Proposition}[Alesker \cite{alesker_survey07}] \label{prop_smooth_vals}
If $G$ is a subgroup of $SO(V)$ acting transitively on the unit sphere of $V$, then the subspace
$\Val^G \subset \Val$ of $G$-invariant elements is a subspace of $\Val^{sm}$.
\end{Proposition}

Valuations which are homogeneous of degree $n-1$ can be explicitly described as follows. 

\begin{Theorem}[McMullen \cite{mcmullen80}] \label{thm_mcmullen_char}
Let $\mu \in \Val_{n-1}$. Then there is a continuous, $1$-homogeneous
function
$f:V^*\rightarrow\R$
such that for all $K \in \mathcal{K}(V)$
\begin{displaymath}
 \mu(K)=V(K[n-1],f).
\end{displaymath}
The function $f$ is unique up to a linear map.
\end{Theorem}


\section{Characterization of the complex projection body}

\begin{Lemma} \label{lemma_invariance_affine}
 Let $W$ be a complex vector space of (complex) dimension $m \geq 2$. Let $G:=\SL(W,\C)$. Let $f:W \to \R$ be
a
continuous function with the property that $f \circ g-f$ is a linear function $l_g$ for each $g \in G$. Then
$f$ is affine. 
\end{Lemma}

\proof
Let us first treat the case $m=2$. Write $f=f^++f^-$ with $f^+$ an even function and
$f^-$ an odd function. Comparing even and odd parts in the equation $f \circ g- f=l_g$ yields that $f^+$ is
$G$-invariant and that $f^- \circ g-f^-=l_g$. Since $G$ acts transitively on $W \setminus \{0\}$, we obtain
that $f^+$ is constant. On the other hand, taking $g=-Id \in \SL(W,\C)$ gives us 
\begin{displaymath}
 f^-=-\frac12 l_{-Id},
\end{displaymath}
i.e. $f^-$ is linear. 

Now let $m >2$. We want to show that $f(x+y)+f(x-y)-2f(x)=0$ for all $x,y \in W$. Let $\tilde W$ be a
two-dimensional complex vector space containing $x$ and $y$ and set $H:=\SL(\tilde W,\C) \cong
\SL(2,\C)$. Each element of
$H$ can be extended to an element of $G$. It follows that the restriction of $f$ to $\tilde W$ satisfies the
assumption of the lemma in the case $m=2$, which we have already discussed. Therefore $f|_{\tilde W}$
is affine, which implies that $f(x+y)+f(x-y)-2f(x)=0$.
\endproof

A valuation $Z:\mathcal{K}(W) \to \mathcal{K}(W^*)$ is said to be {\it $S^1$-bi-invariant} if $Z(qK)=Z(K)$ and $qZ(K)=Z(K)$ for $q\in S^1$ and $K\in\mathcal{K}(W)$.

\begin{Lemma}\label{lemma_homog}
Let $W$ be a complex vector space of complex dimension $m \geq 3$. Let $Z:\mathcal{K}(W) \to \mathcal{K}(W^*)$
be a continuous translation invariant $S^1$-bi-invariant and $\SL(W,\mathbb{C})$-contravariant Minkowski
valuation with degree of homogeneity $k$, $1 \leq k<2m-1$. Then $ZK=\{0\}$, $\forall
K \in \mathcal{K}(W)$.
\end{Lemma}

\proof
First note that the $S^1$-bi-invariance and the homogeneity imply that 
\begin{equation} \label{eq_contravariance_gln}
 Z(gK)=|\det g|^ \frac{k+1}{m} g^{-*}ZK, \quad \forall g \in \GL(W,\mathbb{C}). 
\end{equation}
Indeed, any $g$ can be written as $g=g_0 t q$, where $g_0 \in
\SL(W,\mathbb{C})$, $t \in
\mathbb{R}_{>0}, q \in S^1$. From $\det g=t^m q^m$ we deduce that 
\begin{displaymath}
 Z(gK)=Z(g_0 t qK)=g_0^{-*}t^k Z(qK)=g^{-*} t^{k+1} Z(K)=|\det g|^{\frac{k+1}{m}}  g^{-*} Z(K).
\end{displaymath}

Now we distinguish two cases. 

\noindent {\bf Case $k=m-1$.} 

Let $e_1,\ldots,e_m$ be a complex basis of $W$. Given vectors $w_1,\ldots,w_m \in W$, we can compute the
determinant $\det(w_1,\ldots,w_m) \in \C$ with respect to this basis. 

Let us denote by $[w_1,\dots,w_{m-1}]$ the  parallelotope $[0,w_1]+\dots+[0,w_{m-1}]$. We claim that for each
$w$ and each $m-1$-tuple $\{w_1,\dots,w_{m-1}\}$, we have   
\begin{equation}\label{eq_klain_function}
h(Z[w_1,\dots,w_{m-1}],w)=c |\det(w_1,\dots,w_{m-1},w)|,
\end{equation}
where $c=h(Z[e_1,\dots,e_{m-1}],e_m) \in \R$. 

Since both sides are continuous in $w_1,\ldots,w_{m-1},w$, it is enough to show this equation in the case
where $w_1,\ldots,w_{m-1},w$ are independent over $\C$. 

Define $g\in \GL(W,\C)$ by $gw_1=e_1,\ldots,gw_{m-1}=e_{m-1},gw=e_m$. 

Then, using \eqref{eq_contravariance_gln},
\begin{align*}
h(Z[w_1,\dots,w_{m-1}],w) & = h(g^{-*}Z[w_1,\dots,w_{m-1}],gw)\\
& = h(Zg[w_1,\dots,w_{m-1}],gw)|\det g|^{-1}\\
& = h(Z[e_1,\dots,e_{m-1}],e_m)|\det(w_1,\dots,w_{m-1},w)|.
\end{align*}
This proves \eqref{eq_klain_function}.

Let us fix a Hermitian scalar product on $W$ such that $e_1,\ldots,e_m$ is a Hermitian basis. Then $\SU(W)\cong \SU(m)$ is a subgroup of $\SL(W,\C)$. 

Let $W_0$ be the complex vector space generated by $e_1,\ldots,e_{m-1}$. The stabilizer of $\SU(W)$ at $e_m$
can be identified with the group $\SU(W_0) \cong \SU(m-1)$. 

Define a real-valued valuation $\mu$
on $W_0$ by 
\begin{displaymath}
 \mu(K):=h(ZK,e_m), \quad K \in \mathcal{K}(W_0). 
\end{displaymath}

Clearly $\mu$ is a continuous, translation invariant, even valuation of degree $m-1$. By
\eqref{eq_klain_function}, the Klain function of $\mu$ is $\SU(W_0)$-invariant. Since the Klain function of an
even
valuation describes it uniquely, it follows that $\mu$ is $\SU(W_0)$-invariant. 

The group $\SU(W_0)$ acts
transitively on the unit sphere of $W_0$ (note that this is where our assumption $m \geq 3$ is used). By
Proposition \ref{prop_smooth_vals}, $\mu$ is a smooth valuation,
i.e. $\mu \in
\Val_{m-1}^{sm,+}$. In particular, the Klain function of $\mu$ is a smooth function on $\Gr_{m-1}(W_0)$. 
 
Let $\gamma:\R \to \Gr_{m-1}(W_0)$ be the smooth curve given by 
\begin{displaymath}
 \gamma(t):=\spa_\R \left\{ \cos t e_1 + \sin t i e_2, e_2,e_3,\ldots,e_{m-1}\right\}. 
\end{displaymath}

By \eqref{eq_klain_function} we have 
\begin{displaymath}
 \Kl_\mu(\gamma(t))=c|\det(\cos t e_1 + \sin t i e_2, e_2,e_3,\ldots,e_{m-1},e_m)|=c |\cos t|,
\end{displaymath}
which is smooth only for $c=0$. This implies that the Klain function of the valuation $\mu(K)=h(ZK,w)$ vanishes for each $w\in W$, which by Klain's injectivity theorem implies that $ZK=\{0\}$, i.e. $Z$ is trivial.  

\noindent {\bf Case $1\leq k<m-1$ or $m\leq k<2m-1$.} 

Let $e_1,\ldots,e_m$ be a complex basis of $W$. 

Let $E \subset W$ be the real subspace generated by $e_1,\ldots,e_k$ if $k \leq m$ and by
$e_1,\ldots,e_m, ie_1,\ldots,i e_{k-m}$ if $k>m$.

Let $g \in
\GL(W,\mathbb{C})$ be defined by $ge_j=\lambda_j e_j$ with $\lambda_1,\ldots,\lambda_m \in \R_{>0}$. Let
$D$ be the determinant of the restriction of $g$ to $E$ (considered as an element of $\GL(E,\R)$). Explicitly,
$D=\prod_{j=1}^k \lambda_j$ if $k \leq m$ and $D=\prod_{j=1}^{k-m} \lambda_j^2 \prod_{j=k-m+1}^m \lambda_j$ if
$k>m$. 

Let $w=e_j$ or $w=ie_j$ for some $1 \leq j \leq m$. By Klain's result, the restriction of $h(Z(\cdot),w)$ to
$E$ is a multiple of the $k$-dimensional volume. It follows that 
\begin{equation} \label{eq_klain_fct1}
 h(ZgK,w)=D h(ZK,w), \quad K \in \mathcal{K}(E). 
\end{equation}

On the other hand, 
\begin{equation} \label{eq_klain_fct2}
 h(ZgK,w)=h(ZK,g^{-1}w)  |\det g|^\frac{k+1}{m} =h(ZK,w)  \lambda_j^{-1} |\det g|^\frac{k+1}{m}.
\end{equation}

Comparing \eqref{eq_klain_fct1} and \eqref{eq_klain_fct2}, we get that 
\begin{displaymath}
 D h(ZK,w) \lambda_j =  |\det g|^\frac{k+1}{m} h(ZK,w) 
\end{displaymath}
for all choices of $\lambda_1,\ldots,\lambda_m$. 

The left hand side is clearly polynomial in each $\lambda_j$. Since $|\det g|=\prod
\lambda_j$ and $\frac{k+1}{m} \notin \mathbb{Z}$, the right hand side is a polynomial in the $\lambda_j$ only
if $h(ZK,w)=0$. 

It follows that the support function $h:=h_{ZK}$ vanishes on all lines $\R \cdot e_j, \R \cdot ie_j, j=1,\ldots,m$. Using the convexity of $h$, it follows that $h \equiv 0$ which, by Klain's injectivity theorem, means that $ZK=\{0\}$.
\endproof

\proof[Proof of Theorem \ref{thm_t1}] 

We first check that for each $C \subset \mathbb{C}$, the functional $\Pi_C:\mathcal{K}(W) \to
\mathcal{K}(W^*)$ with
\begin{displaymath}
 h(\Pi_CK,w)=V(K,\ldots,K,C \cdot w), \quad w\in W
\end{displaymath}
satisfies all the stated properties.

It is clear that the function on the right hand side is $1$-homogeneous. For $w_1,w_2 \in W$, we have $C
\cdot (w_1+w_2) \subset C \cdot w_1 + C \cdot w_2$. The monotonicity of the mixed volume implies that 
\begin{displaymath}
 V(K,\ldots,K,C \cdot (w_1+w_2)) \leq V(K,\ldots,K,C \cdot w_1) + V(K,\ldots,K,C \cdot w_2).
\end{displaymath}

Hence the function on the right hand side is the support function of a unique compact convex body $\Pi_CK$
in $W^*$. 
 
Next, we show that $\Pi_C$ is a valuation. By the properties of the mixed volumes, we obtain that
\begin{align*}
 h(\Pi_C(K\cup L) & +\Pi_C(K\cap L),w)  = h(\Pi_C(K\cup L),w)+h(\Pi_C(K\cap L),w)\\
& = V((K \cup L)[2m-1],C \cdot w) +V((K\cap L)[2m-1],C\cdot w)\\
& = V(K[2m-1],C \cdot w)+V(L[2m-1],C \cdot w)\\
& = h(\Pi_CK,w)+h(\Pi_CL,w),
\end{align*}
which implies the valuation property of $\Pi_C$.

The continuity and translation invariance of $\Pi_C$ follow from the corresponding properties for mixed volumes.

To prove the contravariance, we use that mixed volumes are invariant under volume-preserving affine transformations. For each $g
\in \SL(W,\mathbb{C})$ we have
\begin{align*}
h(\Pi_C(gK),w)&=V((gK)[2m-1],C\cdot w)\\
&=V(K[2m-1],g^{-1}C \cdot w)\\
&=V(K[2m-1],C \cdot g^{-1}w)\\
&=h(\Pi_CK,g^{-1}w)\\
&=h(g^{-*}\Pi_CK,w) \quad \text{ by \eqref{eq_adjoint}}.
\end{align*}
It follows that $\Pi_C(gK)=g^ {-*}\Pi_CK$, hence $\Pi_C$ has all the required properties. 

Now let us assume that $Z$ is a continuous translation invariant Minkowski valuation which is
$\SL(W,\mathbb{C})$-contravariant. We want to show that
there exists some compact convex $C \subset \mathbb{C}$ with $Z=\Pi_C$.

We apply the McMullen decomposition \eqref{eq_mcmullen_dec} to $Z$ and write 
\begin{displaymath}
 h(ZK,\cdot)=\sum_{k=0}^{2m}f_k(K,\cdot),
\end{displaymath}
with $f_k(K,\cdot)$ a $1$-homogeneous function. In general it is not known whether $f_k$ is subadditive.
Nevertheless, if $k_0$ and $k_1$ are the minimal and maximal indices $k$ with $f_k \neq 0$, then it is known
(and easy to prove, see \cite{schneider_schuster06}) that $f_{k_0}$ and $f_{k_1}$ are support functions.  

First we show that $k_0=k_1=2m-1$. It is easily checked that degrees $0$ and $2m$ can not
appear since $\Val_0$ is spanned by the Euler characteristic and $\Val_{2m}$ is spanned by the volume. It is
therefore enough to show that there is no non-trivial (i.e. $ZK\neq \{0\}$ for
some $K \in \mathcal{K}(W)$) continuous translation invariant and $\SL(W,\mathbb{C})$-contravariant Minkowski
valuation of
degree $k<2m-1$.

Given such a valuation $Z$, we define
\begin{displaymath}
 \tilde Z(K):=\int_{S^1}\int_{S^1}q_{1}Z(q_{2}K)dq_{1}dq_{2}.
\end{displaymath}

Then $\tilde Z$ is also a non-trivial continuous translation invariant and $\SL(W,\mathbb{C})$-contravariant Minkowski
valuation of degree $k<2m-1$ which is moreover $S^1$-bi-invariant. In
particular, $\tilde Z$ is even, i.e. $\tilde Z(-K)=\tilde Z(K)$ for all $K$. 

By Lemma \ref{lemma_homog} there is no non-trivial continuous translation invariant and
$\SL(W,\mathbb{C})$-contravariant Minkowski valuation $\tilde Z$ of degree $k<2m-1$ which is
$S^1$-bi-invariant. We thus get $k_0=k_1=2m-1$, hence $Z$ must be of degree $2m-1$.

\medskip
McMullen's Theorem \ref{thm_mcmullen_char} implies that for each $w \in W$ there is a continuous
$1$-homogeneous function $f_w:W^* \to \R$
with 
\begin{equation} \label{eq_support_function_zku}
 h(ZK,w)=V(K[2m-1],f_w).
\end{equation}
This function is unique up to a linear function. 

We want to show that $f_w=h_{C\cdot w}$ for some $C \subset \C$ convex, i.e. $f_w(\xi)=h_C((\langle
\xi,w\rangle,\langle\xi, iw\rangle))$ for all $\xi \in W^*$. 

We divide the proof in two steps. In the first step we
show that $f_w(\xi)=G(\langle \xi,w\rangle,\langle\xi, iw\rangle)$ for some $1$-homogeneous function $G:\C \to
\R$.
In the second step we show that $G$ is indeed a support function.

\textsc{Step 1:}
By the contravariance of $Z$, we have for all $g\in \SL(W,\C)$
\begin{displaymath}
 h(ZgK,w)=h(ZK,g^{-1}w)=V(K[2m-1],f_{g^{-1}w})
\end{displaymath}
and 
\begin{displaymath}
 h(ZgK,w)=V((gK)[2m-1],f_w)=V(K[2m-1],f_w \circ g^{-*}) \quad \forall g \in \SL(W,\C).
\end{displaymath}

It follows that 
\begin{equation} \label{eq_contravar_f}
  f_{g^{-1}w} \equiv f_w \circ g^{-*} \quad \forall g \in \SL(W,\C),
\end{equation}
where the equivalence relation means ``up to a linear function''.

Let us fix some non-zero element $w_0 \in W$ and set $f:=f_{w_0}$. Let $H \cong
\SL(m-1,\C)$ be the stabilizer of $\SL(W,\C)$ at $w_0$. Then, by \eqref{eq_contravar_f},  
\begin{equation} \label{eq_inv_up_linear}
 f \circ h^* \equiv f \quad \forall h \in H, 
\end{equation}
i.e. $f$ is $H$-invariant up to linear functions. 

Let $W_0:=\C \cdot w_0$. The inclusion $\iota:W_0 \hookrightarrow W$ induces a projection $\iota^*:W^* \twoheadrightarrow W_0^*$. Note
that each
fiber $(\iota^*)^{-1} \tau, \tau \in W_0^*$ may be identified with the 
complex vector space $W/W_0$ on which $H$ acts. From
\eqref{eq_inv_up_linear} we deduce that $f|_{(\iota^*)^{-1}\tau}$ is $H$-invariant up to linear functions.
By Lemma \ref{lemma_invariance_affine}, $f$ is affine on each fiber of $\iota^*$. 

Since an affine function on an affine subspace not containing the origin is the restriction of some linear
function, and since $f(0)=0$ by homogeneity, it follows that there is a function $\phi:W_0^*
\to W$ such that 
\begin{equation} \label{eq_piecewise_affine}
 f(\xi)=\langle \xi,\phi(\iota^*\xi) \rangle \quad \forall \xi \in W^*. 
\end{equation}

Let $h \in H$. Since $f \circ h^* \equiv f$, there exists some $w_h \in W$ (depending on
$h$)
with 
\begin{displaymath}
 f \circ h^*(\xi)=f(\xi)+\langle \xi,w_h\rangle \quad \forall \xi \in W^*. 
\end{displaymath}
Plugging this into \eqref{eq_piecewise_affine} yields 
\begin{equation} \label{eq_h_phi}
 \langle \xi,h\phi(\iota^* \xi)-\phi(\iota^*\xi)-w_h\rangle=0. 
\end{equation}

Let $\eta \in \ker \iota^*$, i.e. $\iota^*(\xi+\eta)=\iota^*\xi$. Replacing $\xi$ by $\xi+\eta$ in
\eqref{eq_h_phi} yields 
\begin{align*}
 \langle \eta,h\phi(\iota^* \xi)-\phi(\iota^*\xi)-w_h\rangle & =\langle \xi+\eta,h\phi(\iota^*
\xi) -\phi(\iota^*\xi)-w_h\rangle\\
& \quad-\langle \xi,h\phi(\iota^* \xi)-\phi(\iota^*\xi)-w_h\rangle=0.
\end{align*}

Since this is true for all $\eta \in \ker \iota^*$, we get   
\begin{equation} \label{eq_main_invariance}
 h\phi(\iota^* \xi)-\phi(\iota^*\xi)-w_h \in W_0, \quad \forall h \in H, \xi \in W^*.
\end{equation}

Next, we claim that the equivalence class of $\phi(\tau)$ in $W/W_0$ is independent of $\tau \in W_0^*$.  
Let $\xi_1,\xi_2 \in W^*$. Replacing $\xi$ by $\xi_1$ and $\xi_2$ in \eqref{eq_main_invariance} and
subtracting the two equations, we get that the image of $\phi(\iota^* \xi_1)-\phi(\iota^* \xi_2)$ in
$W/W_0$ is fixed by each element $h \in H$. This implies that $\phi(\iota^* \xi_1)-\phi(\iota^* \xi_2) \in
W_0$. The claim thus follows from the surjectivity of $\iota^*$. 

We therefore have some $w' \in W$ with $\phi(\tau)-w' \in W_0$ for all $\tau \in
W_0^*$. Note that the function $\xi \mapsto \langle \xi, \phi(\iota^* \xi)-w'\rangle$ is constant along the
fibers of $\iota^*$.
Therefore, there is some function $p:W_0^* \to \R$ with $\langle \xi, \phi(\iota^* \xi)-w'\rangle=p(\iota^*
\xi)$. Plugging this into \eqref{eq_piecewise_affine} yields 
\begin{displaymath}
 f(\xi)=p(\iota^* \xi)+\langle \xi,w'\rangle. 
\end{displaymath}
In other words, up to a linear function, $f$ is constant along the fibers of
$\iota^*$. Since $f$ was only defined up to a linear function, we may assume from the beginning that $f$ is
constant along the fibers of $\iota^*$, i.e. 
\begin{equation} \label{eq_f_constant_on_fibres}
 f(\xi)=p(\iota^* \xi). 
\end{equation}

Since $\iota^* \xi \in W_0^*$ is determined by its value on two linearly independent vectors, we may rewrite
\eqref{eq_f_constant_on_fibres} as 
\begin{displaymath}
 f(\xi)=G(\langle \xi,w_0 \rangle,\langle \xi,i \cdot w_0\rangle),
\end{displaymath}
where $G:\R^2 \to \R$ is some $1$-homogeneous, continuous function. In particular, $f$ is $H$-invariant. 

Recall that $\SL(W,\C)$ acts
transitively on $W \setminus \{0\}$. Therefore, setting 
\begin{equation} \label{eq_invariance_f}
 f_{gw_0}:=f_{w_0} \circ g^*,\quad \forall g \in \SL(W,\C)
\end{equation}
and $f_0:=0$ yields a well-defined family of continuous $1$-homogeneous functions $f_w:W^* \to \R, w \in W$
with 
\begin{displaymath}
 h(ZK,w)=V(K[2m-1],f_w) \quad \forall w \in W.
\end{displaymath}

More explicitly, if $w=gw_0$, then 
\begin{align*}
 f_w(\xi) & = f_{gw_0}(\xi) \\
& = f_{w_0}(g^* \xi)\\
& = G(\langle g^*\xi,w_0 \rangle,\langle g^*\xi,i \cdot w_0\rangle)\\
& = G(\langle \xi,w \rangle,\langle \xi,i \cdot w\rangle).
\end{align*}

\textsc{Step 2:}
Since $G$ is $1$-homogeneous, in order to see that $G$ is the support function of some compact convex set $C \subset \R^2=\C$, it is enough to show that 
\begin{equation} \label{eq_convexity_g}
 G(z_1+z_2) \leq G(z_1)+G(z_2), \quad \forall z_1,z_2 \in \C. 
\end{equation}

Fixing a complex basis of $W$, we identify $W$ with $\C^m \cong \R^{2m}$. Let $z_3:=-z_1-z_2$ and  
\begin{displaymath}
 w_1:=(\bar z_1,\bar z_2,0,\ldots,0), w_2:=(\bar z_2,\bar z_1,0,\ldots,0), w_3:=(\bar z_3,\bar
z_3,0,\ldots,0).
\end{displaymath}
Note that $w_1+w_2+w_3=0$. 

Let $w_j=r_j \xi_j$ with $r_j \in[0,\infty), \xi_j \in S^{2m-1}$. Define a measure $\rho$ on $S^{2m-1}$ by 
\begin{displaymath}
 \rho:=\sum_{j=1}^3 r_j \delta_{\xi_j}.
\end{displaymath}
Our assumptions imply that  
\begin{displaymath}
 \int_{S^{2m-1}} \xi d\rho(\xi)=\sum_{j=1}^3 r_j \xi_j=\sum_{j=1}^3 w_j=0.  
\end{displaymath}
We approximate $\rho$ weakly by a sequence of measures $\rho_l$ with 
\begin{displaymath}
 \int_{S^{2m-1}} \xi d\rho_l(\xi)=0
\end{displaymath}
which are not concentrated on any great sphere. By Minkowski's existence theorem 
(cf. \cite[Theorem 7.1.2]{schneider_book93}), there exists a convex body $K_l$ with surface
area measure $S_{2m-1}(K_l,\cdot)=\rho_l$. 

For any $u \in W$ we have  
\begin{multline*}
\lim_{l \to \infty} h(ZK_l,u)=\lim_{l \to \infty} \int_{S^{2m-1}} f_u(\xi)dS_{2m-1}(K_l,\xi)=\lim_{l \to
\infty} \int_{S^{2m-1}} f_u(\xi)d\rho_l(\xi)\\
=\int_{S^{2m-1}}
f_u(\xi) d\rho(\xi)=\sum_{j=1}^3 f_u(w_j)=\sum_{j=1}^3 G(\langle u,w_j\rangle,\langle u, iw_j\rangle). 
\end{multline*}

Since $h(ZK_l,\cdot)$ is a support function, we have
\begin{displaymath} 
 h(ZK_l,u_1+u_2) \leq h(ZK_l,u_1)+h(ZK_l,u_2)
\end{displaymath}
for all $u_1,u_2 \in W$. 

Taking the limit yields 
\begin{multline*}
 \sum_{j=1}^3 G(\langle
u_1+u_2,w_j\rangle,\langle u_1+u_2,i w_j\rangle) \leq \sum_{j=1}^3 G(\langle
u_1,w_j\rangle,\langle u_1,i w_j\rangle)\\
+\sum_{j=1}^3 G(\langle
u_2,w_j\rangle,\langle u_2,i w_j\rangle).  
\end{multline*}
In the special case $u_1:=(1,0,\ldots,0)$, $u_2:=(0,1,0,\ldots,0)$, this inequality is \eqref{eq_convexity_g}. 

Hence $G$ is the support function of some convex set $C \subset \C$, i.e.  
\begin{displaymath}
 G(\alpha,\beta)=\sup_{(c_1,c_2) \in C} (\alpha c_1+\beta c_2) 
\end{displaymath}
and therefore
\begin{displaymath}
 f_u(\xi) = G(\langle \xi,u \rangle,\langle \xi,i u\rangle) = \sup_{c \in C} \langle \xi,cu\rangle = h(C
\cdot u,\xi). 
\end{displaymath}
Hence $f_w=h_{C \cdot w}$ and 
\begin{displaymath}
 h(ZK,w)=V(K[2m-1],h_{C \cdot w})=V(K[2m-1],C \cdot w) \quad \forall w \in W, 
\end{displaymath}
which finishes the proof of Theorem \ref{thm_t1}.
\endproof

The assumption $m \geq 3$ is essential in Theorem \ref{thm_t1}. In
the case $m=2$, there are continuous, translation invariant,
$\SL(W,\C)$-contravariant valuations which are not of the form \eqref{eq_def_pi_c}. We do not have a complete
classification in this case, but the following class of examples. 

Fix a basis on $W$. For a compact
convex body $K\subset W$, we denote $\det(K,w):=\{\det(k,w)\,|\,k \in K\}$ which is a
compact convex set in $\mathbb{C}$. 
\begin{Proposition}
\label{prop_n2}
Let $\dim_{\C}W=2$. Let $\mu$ be a continuous, translation invariant, monotone valuation of degree $1$ on
$\C$. Then the operator $Z:\mathcal{K}(W) \to \mathcal{K}(W^*)$ defined by
\begin{displaymath}
 h(ZK,w)=\mu(\det(K,w)),\quad K \in \mathcal{K}(W)
\end{displaymath}
is a continuous, translation invariant, $\SL(W,\C)$-contravariant Minkowski valuation. 
\end{Proposition}

\proof
The map $W \to \C, v \mapsto \det(v,w)$ is linear for each fixed $w \in W$. Hence the image $\det(K,w)$ of $K$
under this map is compact and convex again. Let $\mu$ be a monotone, translation invariant,
continuous valuation of degree $1$. Note that $\mu$ is positive and Minkowski additive (cf. \cite[Theorem 3.2]{goodey_weil84}), i.e. $\mu(K+L)=\mu(K)+\mu(L)$. 

Given $w_1,w_2 \in W$, we have $\det(K,w_1+w_2) \subset \det(K,w_1)+\det(K,w_2)$ and
hence by monotonicity of $\mu$ 
\begin{displaymath}
 \mu(\det(K,w_1+w_2)) \leq \mu(\det(K,w_1)+\det(K,w_2))= \mu(\det(K,w_1))+\mu(\det(K,w_2)).
\end{displaymath}

The function $w \mapsto \mu(\det(K,w))$ is thus the support function of some compact
convex subset $ZK \subset W^*$. 

If $K,L, K \cup L \in \mathcal{K}(W)$, then 
\begin{align*}
\mu(\det(K \cup L,w))+\mu(\det(K \cap L,w)) & =\mu(\det(K,w) \cup \det(L,w)) \\
& \quad + \mu(\det(K,w) \cap
\det(L,w))\\
& = \mu(\det(K,w))+\mu(\det(L,w)).
\end{align*}
It follows that $Z$ is a translation invariant, continuous Minkowski valuation of degree $1$.

To show that $Z$ is $\SL(W,\C)$-contravariant, we compute for $g \in \SL(W,\C)$ and $w \in W$
\begin{displaymath}
 h(ZgK,w)=\mu(\det(gK,w))=\mu(\det(K,g^{-1}w))=h(ZK,g^{-1}w)=h(g^{-*}ZK,w).
\end{displaymath}
\endproof


\section{Geometric inequalities}

Let $V$ be an $n$-dimensional vector space endowed with a volume measure. The classical {\it Brunn-Minkowski
inequality} relates the volume of two convex sets with the volume of its Minkowski sum. If $K,L
\in \mathcal{K}(V)$, and
$0 \leq \lambda \leq 1$, then
\begin{displaymath}
\vol((1-\lambda)K+\lambda L)^{1/n} \geq (1-\lambda)\vol(K)^{1/n}+\lambda\vol(L)^{1/n},
\end{displaymath}
with equality, for some $\lambda \in (0,1)$, if and only if $K$ and $L$ either lie in parallel hyperplanes or
are homothetic. 

We shall use the following generalizations of the Brunn-Minkowski inequality. 

Let $V$ be endowed with a Euclidean scalar product. If $K,L \in \mathcal{K}(V)$ and
$0\leq i\leq n-2$,
then 
\begin{equation}\label{eq_BM3}
W_i(K+L)^{1/(n-i)}\geq W_i(K)^{1/(n-i)}+W_i(L)^{1/(n-i)},
\end{equation}
with equality if and only if $K$ and $L$ are homothetic (cf. \cite{lutwak85}). 

If $K,L,K_1,\dots,K_i \in \mathcal{K}(V)$, $\mathbf{K}=(K_1,\dots,K_i)$ and $0 \leq i \leq n-2$, then 
\begin{equation}\label{BM2}
V((K+L)[n-i],\mathbf{K})^{1/(n-i)} \geq V(K[n-i],\mathbf{K})^{1/(n-i)}+V(L[n-i],\mathbf{K})^{1/(n-i)}. 
\end{equation}
In this inequality the equality conditions are not known.  

A special case of the Aleksandrov-Fenchel inequality is the following
generalized version of the {\it Minkowski inequality}. If $K,L \in \mathcal{K}(V)$ and $0 \leq i \leq n-2$,
then
\begin{equation}\label{eq_M2} 
W_i(K,L)^{n-i} \geq W_i(K)^{n-i-1}W_i(L),
\end{equation}
with equality if and only if $K$ and $L$ are homothetic (cf. \cite[Theorems 6.4.4, 6.6.9]{schneider_book93}). 

Let us return to the complex case. The key lemma for the
proof of the Brunn-Minkowski, Aleksandrov-Fenchel, and Minkowski type inequalities stated in Theorem
\ref{thm_geom_ineq} is the following symmetry property of the mixed
projection body.

\begin{Proposition} \label{prop_symmetry}
Let $K_1,\dots,K_{2m-1} \in \mathcal{K}(W), L_1,\dots,L_{2m-1} \in \mathcal{K}(W^*)$, and $C\subset\C$ a
convex body. Then 
\begin{displaymath}
V(K_1,\dots,K_{2m-1},\Pi_C(L_1,\dots,L_{2m-1}))=V(L_1,\dots,L_{2m-1},\Pi_{\bar C}(K_1,\dots,K_{
2m-1})),
\end{displaymath}
where $\bar C$ is the complex conjugate of $C\subset\C$.
\end{Proposition}

Note that this equation does not depend on the choice of volumes on $W$ and $W^*$. 

\proof
Let us endow $W$ with a Hermitian scalar product and identify $W$ with $W^*$. Set
$\mathbf{K}:=(K_1,\ldots,K_{2m-1})$ and $\mathbf{L}:=(L_1,\ldots,L_{2m-1})$. Then, by the definition of
$\Pi_C$ and by \cite[Theorem 5.1.6]{schneider_book93}
\begin{align*}
V(\mathbf{K},\Pi_C \mathbf{L}) &
=\frac{1}{2m}\int_{S^{2m-1}} h(\Pi_C\mathbf{L},\xi)dS(\mathbf{K},\xi)\\
& = \frac{1}{2m}\int_{S^{2m-1}} \int_{S^{2m-1}} h(C\cdot\xi,w) dS(\mathbf{L},w) dS(\mathbf{K},\xi).
\end{align*}

The statement of the proposition now follows from Fubini's theorem and the relation  
\begin{equation}\label{rel_h_C}
 h(C\cdot \xi,w)=h(\overline C\cdot w,\xi), \quad \xi, w \in W.
\end{equation}
\endproof 

The proof of Theorem \ref{thm_geom_ineq} follows from the symmetry property stated in Proposition
\ref{prop_symmetry} by arguments which were mainly developed by Lutwak
\cite{lutwak86,lutwak_vol_mixed,lutwak93}, see also \cite{schuster07}. Let us only give the proof for
the Brunn-Minkowski type inequality. 

\proof[Proof of Theorem \ref{thm_geom_ineq}.\ref{thm_BMCT})]
Let $Q \in \mathcal{K}(W^*)$. Using Proposition \ref{prop_symmetry}, the general Minkowski
inequality \eqref{eq_M2} and Brunn-Minkowski inequality \eqref{BM2} we have 
\begin{align*}
W_0(Q,\Pi_C(K+L))^\frac{1}{2m-1} & = W_0(K+L,\Pi_{\bar C}Q)^\frac{1}{2m-1} \\
& \geq W_0(K,\Pi_{\bar C}Q)^\frac{1}{2m-1}+W_0(L,\Pi_{\bar C}Q)^\frac{1}{2m-1}\\
& = W_0(Q,\Pi_CK)^\frac{1}{2m-1}+W_0(Q,\Pi_CL)^\frac{1}{2m-1}\\
& \geq W_0(Q)^\frac{1}{2m} \left[W_0(\Pi_CK)^\frac{1}{2m(2m-1)}+W_0(\Pi_CL)^\frac{1}{2m(2m-1)}\right].
\end{align*}
Taking $Q=\Pi_C(K+L)$ we obtain the desired inequality.

The functional $K\mapsto W_{2m-1}(\Pi_{C}K)$ is a translation invariant continous and $U(m)$-invariant valuation of degree $2m-1$. Hence it is a multiple of $W_{1}$ (see \cite{alesker03_un}). Since $C$ is not a point, $\Pi_{C}$ is not trivial and thus there exist a constant $c>0$ such that
\begin{equation} \label{eq_w_equation}
 W_{2m-1}(\Pi_CK)=c W_1(K).
\end{equation}

Now suppose that equality holds in Theorem \ref{thm_geom_ineq}.\ref{thm_BMCT}). Since equality in \eqref{eq_M2} holds only if the two convex bodies are homothetic, the three convex bodies $Q=\Pi_C(K+L), \Pi_CK, \Pi_CL$ are homothetic, i.e. 
\begin{displaymath}
 \Pi_CK=\lambda_1 \Pi_C(K+L)+\xi_1,\quad \Pi_CL=\lambda_2 \Pi_C(K+L)+\xi_2
\end{displaymath}
with $\xi_1,\xi_2 \in W^*$ and $\lambda_1^{2m-1}+\lambda_2^{2m-1}=1$. 
Applying $W_{2m-1}$ to these equations and using \eqref{eq_w_equation} yields
\begin{displaymath}
 W_1(K)=\lambda_1 W_1(K+L),\quad W_1(K)=\lambda_2 W_1(K+L). 
\end{displaymath}
Therefore we have equality in \eqref{eq_BM3}, which implies that $K$ and $L$ are homothetic. 
\endproof


\end{document}